\documentclass[lettersize,journal]{IEEEtran}
\usepackage{amsmath,amsfonts}
\usepackage{algorithmic}
\usepackage{array}
\usepackage[caption=false,font=normalsize,labelfont=sf,textfont=sf]{subfig}
\usepackage{textcomp}
\usepackage{stfloats}
\usepackage{url}
\usepackage{verbatim}
\usepackage{graphicx}
\def\BibTeX{{\rm B\kern-.05em{\sc i\kern-.025em b}\kern-.08em
		T\kern-.1667em\lower.7ex\hbox{E}\kern-.125emX}}
\usepackage{balance}

\RequirePackage[OT1]{fontenc}

\numberwithin{equation}{section}
\newtheorem{thm}{Theorem}

\newtheorem{lem}{Lemma}

\newcommand{\bq}{{\bf q}}
\newcommand{\wht}{\widehat}

\begin{document}

\title{Optimal Sequential Detection by Sparsity Likelihood}
\author{Jingyan Huang}

\maketitle
		
\begin{abstract}
Consider the problem on sequential change-point detection on multiple data streams. We provide the asymptotic lower bounds of the detection delays at all levels of change-point sparsity and we derive a smaller asymptotic lower bound of the detection delays for the case of extreme sparsity. A sparsity likelihood stopping rule based on sparsity likelihood scores is designed to achieve the optimal detections. A numerical study is also performed to show that the sparsity likelihood stopping rule performs well at all levels of sparsity.
We also illustrate its applications on non-normal models.
\end{abstract}


\section{Introduction}

Consider $N$ data streams and $X^{n}_{t}$ denote the observation of the $n$th data stream at time $t$. A change-point $\nu \geq 1$ is said to occur, if for a non-empty subset $\mathcal{N} \subset \{ 1,\cdots, N \} $, the post-change observations $X^{n}_{t}$ for $n \in \mathcal{N}$ (and $t \geq \nu$) have the distributions $f^n_{1}$ which are different from that of pre-change observations, $f^n_{0}$. The problem is to detect the change-point $\nu$ as early as it occurs, while keeping the false alarm rate as low as possible. 

In the literature, sequential change-point detection is defined as a stopping rule $T$ with respect to $\{(X^{1}_{t}, \cdots, X^{N}_{t})\}_{t \geq 1}.$ The interpretation of $T$ is that, when $T=t$, we stop at time $t$ and indicate that an event has occurred in a subset of data streams somewhere in the first $t$ time steps. In the literature, the false alarm rate is usually measured by $1/E_{\infty}T$, where $E_{\infty}T$ is the average run length (ARL) to false alarm when there is no change-point. Our goal is to find a stopping rule $T$ that minimizes the detection delay subject to $E_{\infty}T \geq \gamma,$ where $\gamma$ is a prescribed constant to control the global false alarm rate. Applications for this multi-stream sequential change-point detection (online change-point detection) problem include hospital management, infectious-disease modeling, quality control, surveillance, health care, security and environmental science.

\cite{TV08} considered the setting where sensors detect changes locally and transmit their decisions to a fusion center. In this setting, the authors consider the minimax, uniform and Bayesian formulations for sequential detection in multi-stream data. 
They assume all data streams undergo the change-point if it occurs. Thus, these procedures may include noise from the non-signal data streams when only a subset of data streams undergo change, which may lead to large detection delays. They showed that their procedure is asymptotically optimal with $N$ fixed and $\gamma \rightarrow \infty$.

In analogy to the CUSUM statistic considered by $\cite{Page54, Page55}$ and $\cite{Lorden71}$, $\cite{Mei10}$ proposed a sum of CUSUM stopping rule.
The advantages of his stopping rule are that the distribution changes are not assumed to have occurred simultaneously, and the use of efficient recursive computations. He showed that his stopping rule achieves optimal detection as $\gamma \rightarrow \infty$. However the detection delay is relatively large for sparse signals across sequences.

$\cite{XS13}$ investigated the setting where the proportion of data streams undergoing distribution changes is small. They developed a mixture likelihood ratio (MLR) stopping rule.
$\cite{XS13}$ showed via simulation studies that their MLR stopping rules achieved better detection performance compared to other known approaches over a wide range of $\#\mathcal{N}$. They also provided analytical approximations to ARL of their stopping rules, which are useful for calculating a suitable detection threshold. Analytical approximations to detection delays of their stopping rules were also derived. However an optimality theory for the case of small or moderate $\#\mathcal{N}$ was not provided.

$\cite{Chen19}$ proposes a new framework for sequential change-points detection. 
The approach utilizes nearest neighbor information and can be applied to sequences of multivariate observations
or non-Euclidean data objects. 
One stopping rule is recommended. 
The asymptotic properties of this stopping rule is studied as $L \to \infty,$ where $L$ denotes the number of observations utilized in the stopping rule.
The author also derives an accurate analytic approximation to ARL. Simulations show that the new approach has better performance than likelihood-based approaches for high dimensional data.

In $\cite{Chan17}$, Mei's stopping rule was improved by applying a detectability score transformation on each CUSUM score. 
The advantage of this stopping rule is the suppression of the noise from non-signal data streams. $\cite{Chan17}$ also proposed a modified version of the mixture likelihood ratio (MLR) stopping rule that tests against the limits of detectability. 
The asymptotic lower bounds of the detection delays at all levels of change-point sparsity were provided in $\cite{Chan17}$. For large $\# \mathcal{N}$, the lower bound is trivially given by 1. For moderate $\# \mathcal{N}$, the lower bound grows logarithmically with $N$. For small $\# \mathcal{N}$, the detection delay grows polynomially with $N$. The modified MLR stopping rule achieves optimal detection on all three domains of $\mathcal{N}$.

In this work, a smaller asymptotic lower bound of the detection delays than that in $\cite{Chan17}$ for small $\# \mathcal{N}$ is shown. We also introduce a stopping rule based on sparsity likelihood (SL) score that achieves the detection lower bounds. Another advantage of SL stopping rule is that it is applicable to non-normal models. The calculations of SL scores rely on the calculations of p-values. The p-values could be generated under the models other than normal models.

The outline of this paper is as follows. 
In Section $\ref{model and assumption}$ we introduce the model and assumption on sequential change-point detection in multiple data streams.
In Section $\ref{detection delay lower bound}$ we provide the asymptotic lower bounds of the detection delays at all levels of change-point sparsity, with a smaller asymptotic lower bound than that in $\cite{Chan17}$ at extreme levels of sparsity.
In Section $\ref{optimal detection using sparsity likelihood score}$ we show that a SL stopping rule achieves optimal detection at all levels of change-point sparsity. A window-limited rule proposed by $\cite{Lai95}$ integrates into the stopping rule to show computational efficiency.
In Section $\ref{numerical study}$ a simulation study is performed to show that SL stopping rules and the modified MLR stopping rules are relatively competitive over the full range of $\mathcal{N}$. Furthermore, SL stopping rule outperforms other competitors for small $\# \mathcal{N}$. In section 6, we discuss the extensions of SL stopping rule to other non-normal models. 
In the Appendix, proofs are provided to show that SL stopping rules are optimal.

\section{Model and assumption}  \label{model and assumption}

Consider $N$ data streams. Let $X^{n}_{t}$ denote $t$th observation of the $n$th data stream. We assume first that $X^{n}_{t}$ are independent normal with unit variance. Let $\mu^{n}_{t}$ denote the mean of $X^{n}_{t}$. Hence, 
\begin{equation}  \label{normod}
	X_{t}^{n} \sim N(\mu^{n}_{t},1).
\end{equation}
Assume that at some unknown time, $\nu \geq 1,$ there are mean shifts in a subset $\mathcal{N}$ of the data stream. More specifically, we assume that
\begin{equation}  \label{munor}
	\mu^{n}_{t} = \Delta \mathbf{1}_{\{t \geq \nu, n \in \mathcal{N}\}}, \mbox{ for some } \Delta \neq 0,
\end{equation}
with $\mathbf{1}_{\{1 \in \mathcal{N}\}}, \mathbf{1}_{\{2 \in \mathcal{N}\}}, \cdots, \mathbf{1}_{\{N \in \mathcal{N}\}}$ i.i.d. Bernoulli($\epsilon$) for some $0 < \epsilon < 1$. We shall let $P_{\nu} (E_{\nu})$ denote probability measure (expectation) with respect to distribution changes at time $\nu$($E_{\nu}$ denotes expectation with respect to $\# \mathcal{N} \sim $Binomial$(N,\epsilon)$), with $\nu = \infty$ indicating no changes.

A standard measure of the performance of a stopping rule $T$ is the (expected) detection delay:
\begin{equation}  \label{delay}
	D_{N}(T) = \sup_{1 \leq \nu < \infty} E_{\nu}(T-\nu+1|T \geq \nu), 
\end{equation}
subject to the constraint that $ARL(T)(:= E_{\infty}T) \geq \gamma$ for some $\gamma \geq 1$. Often but not always the worst case detection delay of a stopping rule occurs at $\nu=1$, i.e. it is more difficult to detect when a change occurs at earlier stages rather than at later stages.

Let $E_{\nu, \mathcal{N}}$ denote expectation with respect to $X_{t}^{n} \sim N(\mu_{t}^{n},1),$ with $\mu^{n}_{t} = \Delta \mathbf{1}_{\{t \geq \nu, n \in \mathcal{N}\}} \mbox{ for some } \Delta \neq 0.$ That is, $E_{\nu,\mathcal{N}}$ denotes expectation with respect to a fixed $\mathcal{N}$. For a given stopping rule $T$, a standard measure of the performance of a stopping rule $T$ is the (expected)detection delay which is defined as
\begin{equation}
	D_{N}^{(V)}(T) = \sup_{1 \leq \nu < \infty} \Big[\max_{\mathcal{N}:\# \mathcal{N}=V } E_{\nu, \mathcal{N}}(T-\nu+1|T \geq \nu) \Big],
\end{equation}
where the average run length (ARL) to false alarm when there is no change-point  is subject to the constraint that ARL$(T)(:= E_{\infty}T) \geq \gamma$ for some $\gamma \geq 1$. 

\bigskip{}
\section{Detection delay lower bound} \label{detection delay lower bound}

We show here the (asymptotic) lower bounds of detection delay, under the conditions that as $N \rightarrow \infty,$
\begin{equation}  \label{gam}
	\log \gamma \sim N^{\zeta} \textrm{ for some }   \textrm{$0 < \zeta \leq 1$},
\end{equation}
\begin{equation} \label{ep}
	\epsilon \sim N^{-\beta} \textrm{ for some }   \textrm{$0 < \beta < 1$},
\end{equation}
where $a_{n} \sim b_{n}$ if $\lim\limits_{n \rightarrow \infty} (a_{n}/b_{n})=1.$

For $\beta < \frac{1-\zeta}{2},$ the sparsity likelihood stopping rules achieve asymptotic detection delay of 1, and are hence optimal.

It was shown in $\cite{Chan17}$ that for $\frac{1-\zeta}{2} < \beta < 1-\zeta,$ the detection delay lower bound grows logarithmically with $N$. The proportionality constant is 
\begin{eqnarray} \label{rhozb}
	\rho_{Z}(\beta,\zeta) = \left\{ \begin{array}{ll} \beta - \frac{1-\zeta}{2}, & \mbox{if } \frac{1-\zeta}{2}< \beta \leq \frac{3(1-\zeta)}{4}, \\
		(\sqrt{1-\zeta}-\sqrt{1-\zeta-\beta})^{2}, & \mbox{if } \frac{3(1-\zeta)}{4}< \beta \leq 1-\zeta. \end{array} \right.
\end{eqnarray}
This constant, first defined in $\cite{CW15}$, it is a two-dimensional extension of the Donoho–Ingster–Jin constants $\rho_{Z}(\beta):= \rho_{Z}(\beta,0),$ for sparse normal mixture detection, see $\cite{DJ04}$ and $\cite{Ingster97, Ingster98}$. The extension results from the multiple comparisons in detecting a normal mean shift, here for sequential change-point detection on multiple data streams, and in $\cite{CW15}$ for fixed-sample change-point detection on multiple sequences.

For $\beta > 1-\zeta,$ we provide a smaller (asymptotic) lower bound of detection delay compared with that in $\cite{Chan17}$. The results below apply to all stopping rules.

\bigskip
\begin{thm} \label{thm21}
	Assume (\ref{normod}), (\ref{munor}) and (\ref{ep}). Let $T$ be a stopping rule such that $ARL(T) \geq \gamma,$ with $\gamma$ satisfying (\ref{gam}). 
	
	{\rm (a)} If $\frac{1-\zeta}{2} < \beta < 1-\zeta,$ then
	\begin{equation}  \label{Thm1a}
		\liminf_{N \rightarrow \infty} \frac{D_{N}(T)}{\log N} \geq 2\Delta^{-2}\rho_{Z}(\beta, \zeta).
	\end{equation}
	
	{\rm (b)} When $V=o(\frac{N^{\zeta}}{\log N})$ for $\beta > 1-\zeta$, then
	\begin{equation}   \label{Thm1b}
		\liminf_{N \rightarrow \infty} \frac{D_{N}^{(V)}(T)}{N^{\zeta}} \geq 2\Delta^{-2}V^{-1}.
	\end{equation}
\end{thm}

\bigskip

\bigskip{}
\section{Optimal detection using sparsity likelihood score}  \label{optimal detection using sparsity likelihood score}

We consider here the test of $H_{0}: \nu=\infty$ versus $H_{1}: \nu = s$ for some $s \leq t$.\\
Let $k = t-s+1,$ $S^{n}_{st}=\sum\limits^{t}_{i=s}X^{n}_{i}, Z^{n}_{st}=\frac{S^{n}_{st}}{\sqrt{k}},$ and  $p^{n}_{st}=2\Phi(-|Z^{n}_{st}|)$, with $\Phi$ denoting the distribution function of the standard normal. \\
As in $\cite{Hu2022}$, let $f_{1}(p) = \frac{1}{p(2-\log p)^{2}} - \frac{1}{2}$ and $ f_{2}(p) = \frac{1}{\sqrt{p}}-2$. Note that $\int_{0}^{1}f_{i}(p)dp=0$ for $i=1$ and $2.$ Let $\lambda_{1} \geq 0$ and $\lambda_{2} > 0.$ Define the sparsity likelihood score
\begin{eqnarray}   \label{lP}
	\ell(\mathbf{p}) &=& \sum_{n=1}^{N} l(p^{n}), \\ \nonumber
	\mbox{where } \ell(p) &=& \log \Big(1+\frac{\lambda_{1}\log N}{N}f_{1}(p)+\frac{\lambda_{2}}{\sqrt{N\log N}}f_{2}(p)\Big).
\end{eqnarray}
When working with p-values small compared to $\frac{1}{N \log N},$ the selection of $\lambda_{1}>0$ is advantageous because $\frac{\log N}{N}f_{1}(p)$ dominates $\frac{1}{\sqrt{N \log N}}f_{2}(p)$ for $p$ small compared to $\frac{1}{N \log N}.$
We combine the p-values using the score $\ell(\mathbf{p}_{st}),$ where $\mathbf{p}_{st}=(p^{1}_{st}, \cdots, p^{N}_{st})$ and sparsity likelihood stopping rule can be expressed as

\begin{equation}  \label{SL}
	T_{SL} = \inf \Big\{t: \max_{k=t-s+1 \in \mathcal{K}} \ell(\mathbf{p}_{st}) \geq C_{\gamma} \Big\}.  
\end{equation}
We follow the window sizes considered in window-limited GLR stopping rule in $\cite{Lai95}$ for saving computations. Consider
\begin{equation} \label{win}
	\mathcal{K} = \{1, \cdots, k_{1}\}\cup\{ \left \lfloor{r^{j}k_{1}}\right \rfloor: j \geq 1\}, k_{1} \geq 1, r>1.
\end{equation}
In Theorem $\ref{thm22}$ below we show that sparsity likelihood stopping rule achieves minimum detection delay. We select parameters 
\begin{equation}  \label{lambda}
	\lambda_{1} > 0 \mbox{ and } \lambda_{2} = \sqrt{\frac{\log \gamma}{\log\log \gamma}}.
\end{equation}

\bigskip
\begin{thm} \label{thm22} 
	Assume (\ref{normod}) and (\ref{munor}). Consider stopping rule $T_{SL}$, with window sizes in (\ref{win}). If $ARL(T_{SL})=\gamma,$ then the stopping rule threshold $C_{\gamma} \leq  \log (4\gamma^{2}+2\gamma).$ In addition, if (\ref{gam}),(\ref{ep}) hold, then the following hold as $N \rightarrow \infty:$
	
	\medskip
	{\rm (a)} If $\beta<\frac{1-\zeta}{2},$ then
	\begin{equation}
		D_{N}(T_{SL})\rightarrow 1.
	\end{equation}
	
	{\rm (b)} If $\frac{1-\zeta}{2} < \beta < 1-\zeta,$ then
	\begin{equation}
		\frac{D_{N}(T_{SL})}{\log N} \rightarrow 2\Delta^{-2}\rho_{Z}(\beta,\zeta).
	\end{equation}
	
	{\rm (c)} If $V=o(\frac{N^{\zeta}}{ \log N})$ for $\beta > 1-\zeta$, then
	\begin{equation}
		\frac{D_{N}^{(V)}(T_{SL})}{N^{\zeta}} \rightarrow 2\Delta^{-2}V^{-1}.
	\end{equation}
\end{thm}

\bigskip
\section{Simulation}  \label{numerical study}

In this section, we compare the SL stopping rule against the four stopping rules below.

For testing the hypothesis that change-points $\nu=s$ for some $s \leq t,$ the most powerful test at time $t$ is the log likelihood
\begin{equation*}
	\ell_{st}=\sum_{n=1}^{N}\ell^{n}_{st}, \mbox{ where } \ell^{n}_{st}=\log(1-\epsilon+\epsilon e^{\Delta S^{n}_{st}-k \Delta^{2}/2}),
\end{equation*}
with $k=t-s+1$ and $S^{n}_{st}=\sum_{i=s}^{t}X^{n}_{i}$. 

$\cite{XS13}$ considered the situation that $\epsilon$ and $\Delta$ are known. Since the change-point $\nu$ is unknown, $\cite{XS13}$ suggest to maximize $\ell_{st}$ over $s$. To sum up, the stopping rule can be expressed as
\begin{equation}  \label{XS}
	T_{XS}(\epsilon_{0})=\inf \Big\{t:\max_{k=t-s+1 \in \mathcal{K}} \wht{\ell}_{st}(\epsilon_{0}) \geq C_{\gamma} \Big\},
\end{equation}
where $\wht {\ell}_{st}(\epsilon_{0})=\sum\limits_{n=1}^{N} \wht {\ell}^{n}_{st}(\epsilon_{0})$,
\begin{equation*}
	\wht {\ell}^{n}_{st}(\epsilon_{0}) = \log(1-\epsilon_{0}+\epsilon_{0}e^{(Z^{n+}_{st})^{2}/2}), Z^{n}_{st}=S^{n}_{st}/\sqrt{k},
\end{equation*}
and $x^{+}$ denotes the positive part of $x$. Here thresholding by the positive
part plays the role in limiting the current considerations only to the data streams that appear to be affected by the change-point. The set $\mathcal{K}$ in $(\ref{XS})$ refers to a pre-determined set of window sizes. 

$\cite{Mei10}$ introduced the stopping rule:
\begin{equation}  \label{Mei2010}
	T_{Mei}=\inf \Big\{t:\sum_{n=1}^{N}\max_{0<s \leq t}(\Delta_{0}S^{n}_{st}-k\Delta_{0}^{2}/2)^{+} \geq C_{\gamma} \Big\},
\end{equation}
with a pre-determined $\Delta_{0}.$
In $\cite{Chan17}$, Mei's stopping rule was improved by applying a detectability score transformation on each CUSUM score. The advantages of $(\ref{Mei2010})$ are the efficient recursive computation of the stopping rule and the suppression of the noise from non-signal data streams. However, there is information loss. 

Let $R^{n}_{t}$ be the CUSUM score of the $n$th detector at time $t$,
satisfying
\begin{equation}   \label{3.1} 
	R^{n}_{0}=0, R^{n}_{t} = (R^{n}_{t-1}+\Delta_{0}X^{n}_{t}-\Delta_{0}^{2}/2)^{+}, t \geq 1.
\end{equation}
Define
\begin{equation}   \label{3.2} 
	T_{Mei}(\epsilon_{0})=\inf \Big\{t:\sum_{n=1}^{N}g_{M}(R^{n}_{t}) \geq C_{\gamma} \Big\},
\end{equation}
with the detectability score transformation
\begin{equation}   \label{3.3} 
	g_{M}(x)=\log[1+\epsilon_{0}(\lambda_{M}e^{x/2}-1)], \lambda_{M}>0.
\end{equation}
This is an extension of Mei's test, for $T_{Mei}(1)$ is equivalent to $T_{Mei}.$ Define 
\begin{equation}   \label{3.4}
	D_{N,k}(T)=\sup_{k \leq \nu < \infty}E_{\nu}(T-\nu+1|T \geq \nu).
\end{equation}

Motivated by the MLR stopping rule by $\cite{XS13}$, \cite{Chan17} introduces a modified version of MLR stopping rule based on the limits of detectability:
\begin{equation}
	T_{S}(\epsilon_{0}) = \inf \Big\{t: \max_{k=t-s+1 \in \mathcal{K}} \sum_{n=1}^{N}g(Z^{n+}_{st}) \geq C_{\gamma} \Big\}, 
\end{equation}
where $g(z)=\log[1+\epsilon_{0}(\lambda e^{z^{2}/4}-1)], \lambda = 2(\sqrt{2}-1),$ and $\epsilon_{0}=N^{-1/2}, [(\log \gamma)/N]^{1/2}$ which are optimal choices in $\cite{Chan17}$.

\subsection{Assumption and revised SL stopping rule}

The above four stopping rules are designed to test the change of mean $\Delta > 0$. To benchmark the SL stopping rule against the above four stopping rules, we restrict our assumption $\Delta \neq 0$ to $\Delta > 0 $ in $(\ref{munor})$ in Section \ref{numerical study}. More specifically, we assume 
\begin{equation}  \label{munor1}
	\mu^{n}_{t} = \Delta \mathbf{1}_{\{t \geq \nu, n \in \mathcal{N}\}}, \mbox{ for some } \Delta > 0,
\end{equation}
where $\nu \geq 1 $ denote the unknown time when there are mean shifts in a subset $\mathcal{N}$ of the data stream.

For the revised SL stopping rule, we consider the SL score based on the one-sided p-value. Let $k = t-s+1, S^{n}_{st}=\sum\limits^{t}_{i=s}X^{n}_{i}, Z^{n}_{st}=\frac{S^{n}_{st}}{\sqrt{k}},$ and  $p^{n}_{st}=\Phi(-Z^{n}_{st})$. Define the sparsity likelihood score
\begin{eqnarray}   
	\ell(\mathbf{p}) &=& \sum_{n=1}^{N} l(p^{n}), \\ \nonumber
	\mbox{where } \ell(p) &=& \log \Big(1+\frac{\lambda_{1}\log N}{N}f_{1}(p)+\frac{\lambda_{2}}{\sqrt{N\log N}}f_{2}(p)\Big),
\end{eqnarray}
where $f_{1}(p), f_{2}(p), \lambda_{1}, \lambda_{2}$ are defined as in Section $\ref{optimal detection using sparsity likelihood score}$. The p-values are combined using the score $\ell(\mathbf{p}_{st}),$ where $\mathbf{p}_{st}=(p^{1}_{st}, \cdots, p^{N}_{st})$ and sparsity likelihood stopping rule can be expressed as
\begin{equation}  \label{SL}
	T_{SL}^{1} = \inf \Big\{t: \max_{k=t-s+1 \in \mathcal{K}} \ell(\mathbf{p}_{st}) \geq C_{\gamma} \Big\}, 
\end{equation}
and the window sizes considered are in (\ref{win}).

\subsection{Simulation results}

To compare $T_{SL}^{1}$ with other four stopping rules, we firstly provide the optimal choices of parameters for $T_{SL}^{1}$.  For $T_{SL}^{1}$, we select parameters $\lambda_{1}=1$ as the same choice of $\lambda_{1}$ in simulation studies in $\cite{Hu2022}$. $\lambda_{2}$ we choose is $\lambda_{2}=\{0.2,0.4,0.6,\cdots,1.8,\sqrt{\frac{\log \gamma}{\log\log \gamma} } \}$.
	
The first simulation model we consider is the model in $\cite{Chan17}$, with $(\ref{normod})$ and $(\ref{munor1})$. The parameters selected are $N = 100, \Delta=1$ and $\# \mathcal{N}$ ranging from $1$ to $100$. The thresholds $C_{\gamma}$ is calibrated to ARL of 5000. The set of window sizes chosen is $\mathcal{K} = \{1,\cdots,200\}$.

We conduct 500 Monte Carlo trials for the estimation of each average run length. We also conduct 500 Monte Carlo trials to simulate the detection delays when the change-point occurs at time $\nu=1$, i.e. providing the estimated values of the worst case detection delays. 

As Table $\ref{table:22}$ shows, SL stopping rule with $\lambda_{2}=\sqrt{\frac{\log \gamma}{\log\log \gamma}}$ has the smallest detection delays compared to all other parameter settings over \\ $\# \mathcal{N} = \{5, 10, 30, 50, 100\}$. Compared with $\lambda_{2}=\sqrt{\frac{\log \gamma}{\log\log \gamma}}$, the detection delay with $\lambda_{2}=1.0$ decreases obviously for $\# \mathcal{N}=1$ and with the slight effects for $\# \mathcal{N} = \{5,10,30,50,100\}$.
		
We select SL stopping rules with $\lambda_{1} = 1, \lambda_{2}=\sqrt{\frac{\log \gamma}{\log\log \gamma}}$ and $\lambda_{1} = 1, \lambda_{2}=1.0$ to compare against the above four stopping rules.
The thresholds for the stopping rules are in Table $\ref{table:23}$, and the detection delays are in Table $\ref{table:24}$. The simulation outcomes for Xie and Siegmund’s stopping rule, Mei's stopping rules, the extension of Mei's stopping rules and the modified MLR stopping rule are reproduced from $\cite{Chan17}$. 

We see that SL stopping rules and the modified MLR stopping rules have smaller detection delays compared to other competitors over the full range of $\mathcal{N}$. To compare the modified MLR stopping rules with SL stopping rules, we construct two groups. We construct the group that works better for small $\#\mathcal{N}$, consisting of the modified MLR stopping rule with $\epsilon_{0}=N^{-1/2}$ and the SL stopping rule with $\lambda_{1} = 1, \lambda_{2} = 1.0$. We also construct the group that works better for moderate and large $\#\mathcal{N},$ consisting of the modified MLR stopping rule with $\epsilon_{0}=[(\log \gamma)/N]^{1/2}$ and the SL stopping rule with $\lambda_{1} = 1, \lambda_{2} = \sqrt{\frac{\log \gamma}{\log\log \gamma}}$. Table $\ref{table:24}$ shows that for moderate and large $\#\mathcal{N},$ there is not much differences on detection delays between SL stopping rules and the modified MLR stopping rule for each group. For small $\#\mathcal{N},$ that is, $\#\mathcal{N} = 1$ and $3$, SL stopping rule has the smaller detection delays compared to that of MLR stopping rule for each group.

\begin{table}[ht]  
	\centering
	\caption{Thresholds $C_{\gamma}$ for SL stopping rules calibrated to $ARL=5000$.}
	\begin{tabular}{l|ll}
		\hline
		& \multicolumn{2}{c}{$N=100$} \\
		\hline
		$Parameters$ & $C_{\gamma}$ & $ARL$ \\
		\hline
		$\lambda_{1}=1,\lambda_{2}=0.2$ & 6.400  & 5071 \\
		
		$\lambda_{1}=1,\lambda_{2}=0.4$ & 6.430  & 5069 \\
		
		$\lambda_{1}=1,\lambda_{2}=0.6$ & 6.475 & 4951 \\
		
		$\lambda_{1}=1,\lambda_{2}=0.8$ & 6.560 & 4943 \\
		
		$\lambda_{1}=1,\lambda_{2}=1.0$ & 6.650  & 5088\\
		
		$\lambda_{1}=1,\lambda_{2}=1.2$ & 6.760  & 4991 \\
		
		$\lambda_{1}=1,\lambda_{2}=1.4$ & 6.860 & 5024 \\
		
		$\lambda_{1}=1,\lambda_{2}=1.6$ & 6.960  & 5003 \\
		
		$\lambda_{1}=1,\lambda_{2}=1.8$ & 7.060 & 5018 \\
		
		$\lambda_{1}=1,\lambda_{2}=\sqrt{\frac{\log \gamma}{\log\log \gamma}} = 1.99$ & 7.160 & 5036 \\
		\hline
	\end{tabular}
	\label{table:21}
\end{table}

\begin{table}[ht]   
	\centering
	\caption{Detection delays when $\# \mathcal{N}$ (out of $N = 100$) data streams undergo normal distribution changes for SL stopping rule.}
	\begin{tabular}{l|lllllll}
		\hline
		& \multicolumn{7}{|c}{$\# \mathcal{N}$}\\
		\hline
		$Parameters$ & 1 & 3 & 5 & 10 & 30 & 50 & 100\\
		\hline
		$\lambda_{1}=1,\lambda_{2}=0.2$ & 24.5 & 13.5 & 10.4 & 7.1 & 3.8 & 2.8 & 1.8 \\
		$\lambda_{1}=1,\lambda_{2}=0.4$ & 24.7 & 13.3 & 10.1 & 6.7 & 3.4 & 2.4 & 1.4 \\
		$\lambda_{1}=1,\lambda_{2}=0.6$ & 24.8 & 13.2 & 9.9 & 6.4 & 3.1 & 2.2 & 1.2 \\
		$\lambda_{1}=1,\lambda_{2}=0.8$ & 25.3 & 13.3 & 9.7 & 6.2 & 2.9 & 2.0 & 1.1 \\		
		$\lambda_{1}=1,\lambda_{2}=1.0$ & 25.9 & 13.3 & 9.7 & 6.0 & 2.7 & 1.8 & 1.0 \\		
		$\lambda_{1}=1,\lambda_{2}=1.2$ & 26.4 & 13.3 & 9.6 & 5.9 & 2.6 & 1.7 & 1.0 \\
		$\lambda_{1}=1,\lambda_{2}=1.4$ & 26.8 & 13.4 & 9.6 & 5.8 & 2.5 & 1.7 & 1.0\\		
		$\lambda_{1}=1,\lambda_{2}=1.6$ & 27.4 & 13.5 & 9.6 & 5.7 & 2.4 & 1.6 & 1.0\\		
		$\lambda_{1}=1,\lambda_{2}=1.8$ & 28.0 & 13.6 & 9.6 & 5.7 & 2.3 & 1.5 & 1.0\\
		$\lambda_{1}=1,\lambda_{2}=\sqrt{\frac{\log \gamma}{\log\log \gamma}}=1.99$ & 28.6 & 13.7 & 9.6 & 5.6 & 2.2 & 1.5 & 1.0\\
		\hline
	\end{tabular}
	\label{table:22}
\end{table}


\begin{table}[ht]  
	\centering
	\caption{Thresholds $C_{\gamma}$ for stopping rules calibrated to $ARL=5000$. The upper bounds of the thresholds, as given in the statements of Theorems in $\cite{Chan17}$ and in this study, are in brackets.}
	\begin{tabular}{l|ll}
		\hline
		& \multicolumn{2}{c}{$N=100$}\\
		\hline
		$Test$ & $C_{\gamma}$ & $ARL$\\
		\hline
		$Mei$ & 88.500 (106.8) & 4997\\
		$Mei(N^{-1/2})$ & 3.480 (9.81) & 4994 \\
		$Mei(3N^{-1/2})$ & 5.020 (9.61)& 4976 \\
		$S(N^{-1/2})$ & 4.250 (18.42) & 5066\\
		$S(3N^{-1/2})$ & 6.300 (18.42) & 5195 \\
		$SL(\lambda_{1}=1,\lambda_{2}=1.0)$ & 6.650 (18.42)  & 5088\\
		$SL(\lambda_{1}=1,\lambda_{2}=\sqrt{\frac{\log \gamma}{\log\log \gamma}} = 1.99)$ & 7.160 (18.42) & 5036  \\
		\hline
	\end{tabular}
	\label{table:23}
\end{table}

\begin{table}[ht]   
	\centering
	\caption{Detection delays when $\# \mathcal{N}$ (out of $N = 100$) data streams undergo normal distribution changes.}
	\begin{tabular}{l|lllllll}
		\hline
		& \multicolumn{7}{|c}{$\# \mathcal{N}$}\\
		\hline
		$Test$ & 1 & 3 & 5 & 10 & 30 & 50 & 100\\
		\hline
		$XS(1)$ & 52.3 & 18.7 & 12.2 & 6.7 & 3.0 & 2.3 & 2.0\\
		$XS(0.1)$ & 31.6 & 14.2 & 10.4 & 6.7 & 3.5 & 2.8 & 2.0\\
		$Mei$ & 53.2 & 23.0 & 15.7 & 9.6 & 4.9 & 3.8 & 3.0\\
		$Mei(0.1)$ & 26.4 & 14.6 & 10.8 & 7.7 & 4.5 & 3.4 & 2.3\\
		$Mei(0.3)$ & 34.3 & 15.9 & 11.8 & 7.6 & 4.1 & 3.1 & 2.0\\
        $S(0.1)$ & 26.8 & 13.4 & 9.6 & 6.4 & 2.8  & 2.0 & 1.1\\
        $S(0.3)$ & 32.6 & 14.0 & 9.5 & 5.6 & 2.3  & 1.5 & 1.0\\
		$SL(\lambda_{1}=1,\lambda_{2}=1.0)$ & 25.9 & 13.3 & 9.7 & 6.0 & 2.7 & 1.8 & 1.0 \\
		$SL(\lambda_{1}=1,\lambda_{2}=\sqrt{\frac{\log \gamma}{\log\log \gamma}}=1.99)$ & 28.6 & 13.7 & 9.6 & 5.6 & 2.2 & 1.5 & 1.0\\
		\hline
	\end{tabular}
	\label{table:24}
\end{table}


\clearpage

\section{Extensions to non-normal model} 

Here we extend our current models to other non-normal models and apply the SL stopping rule.

The calculations of SL scores rely on the calculations of p-values. The p-values could be generated under the models other than normal models. In this section, we aim to extend normal models to other non-normal models with the application of SL stopping rule.

We consider $N$ data streams. Suppose that $X^{n}_{t}$ is $t$th observation of the $n$th data stream. We assume that $X^{n}_{t}$ are independent random variables from a certain distribution. Let $E(X^{n}_{t}) = \mu^{n}_{t},$ where
\begin{eqnarray*}
	\mu_{t}^{n} & = & \left\{ \begin{array}{ll}  \Delta_{1}, & t \geq \nu, n \in \mathcal{N},\\
		\Delta_{0}, & \mbox{otherwises}, 
	\end{array} \right.
\end{eqnarray*}
where $\Delta_{0}$ is assumed known, $\Delta_{1} \neq \Delta_{0,}$and $\nu$ is the unknown time when the mean shifts in a subset $\mathcal{N}$ of the data stream occur. 

We consider the test of $H_{0}: \nu=\infty$ versus $H_{1}: \nu = s$ for some $s \leq t$. Let $k = t-s+1,$ and $p^{n}_{st}$ can be calculated from the observations $\{(X^{n}_{i})\}_{s\leq i \leq t}$. If $X^{n}_{t}$ follow discrete distribution, there are continuity adjustments for $p^{n}_{st}$ so that it is distributed as Uniform(0,1) under the null hypothesis. $f_{1}(p)$ and $f_{2}(p)$ are defined in Section $\ref{optimal detection using sparsity likelihood score}$. Let $\lambda_{1} \geq 0$ and $\lambda_{2} > 0.$ Define the sparsity likelihood score
\begin{eqnarray}   \label{lP}
	\ell(\mathbf{p}) &=& \sum_{n=1}^{N} l(p^{n}), \\ \nonumber
	\mbox{where } \ell(p) &=& \log \Big(1+\frac{\lambda_{1}\log N}{N}f_{1}(p)+\frac{\lambda_{2}}{\sqrt{N\log N}}f_{2}(p)\Big).
\end{eqnarray}
And we consider the SL stopping rule
\begin{equation}  
	T_{SL} = \inf \Big\{t: \max_{k=t-s+1 \in \mathcal{K}} l(\mathbf{p}_{st}) \geq C_{\gamma} \Big\},
\end{equation}
where $\mathcal{K}$ is the set of window sizes considered in $(\ref{win})$.

In the next two subsections, we consider the application of SL stopping rules on two non-normal models, Poisson model and binomial model.

\subsection{Poisson model and an example} \label{Poi_eg} 

Let $X^{n}_{t}$ denote $t$th observation of the $n$th sequence. We assume that $X^{n}_{t}$ are independent Poisson random variables. Let $\mu^{n}_{t}$ denote mean of $X^{n}_{t}$. Hence,
\begin{eqnarray}   \label{poimod}
	X^{n}_{t} & \sim & \mbox{Poisson}(\mu^{n}_{t}),
\end{eqnarray}
\begin{eqnarray*}
	\mu_{t}^{n} & = & \left\{ \begin{array}{ll}  \Delta_{1} \neq 0.015, & t \geq \nu, n \in \mathcal{N},\\
		\Delta_{0} = 0.015, & \mbox{otherwises}. \end{array} \right.
\end{eqnarray*}

To test $H_{0}: \nu=\infty$ versus $H_{1}: \nu = s$ for some $s \leq t$. 
Let $k = t-s+1,$ $S^{n}_{st}=\sum\limits^{t}_{i=s}X^{n}_{i} \sim \mbox{Poisson}(k\mu^{n}_{t} )$. P-values are calculated based on $p^{n}_{st} = 2 \min (\varphi^{n}_{st}, 1-\varphi^{n}_{st})$, where $\varphi^{n}_{st} \sim$ Uniform$(P( X < S^{n}_{st}),P(X \leq S^{n}_{st}))$, and $X \sim \mbox{Poisson}(k\Delta_{0} )$.

We select $\Delta_{0} = 0.015, \Delta_{1}= 0.3, \lambda_{1}=1, \lambda_{2}= \sqrt{\frac{\log \gamma}{\log \log \gamma}} $. Let $k = t-s+1$ and $S^{n}_{st}=\sum\limits^{t}_{i=s}X^{n}_{i} \sim \mbox{Poisson}(k\Delta_{1})$ if $n \in \mathcal{N}$ and $\nu = 1$. p-values are computed and SL stopping rule is applied. 

\subsection{Binomial model and an example} \label{Bin_eg} 

Let $X^{n}_{t}$ denote $t$th observation of the $n$th sequence. We assume that $X^{n}_{t}$ are independent binomial random variables. Let $p^{n}_{t}$ denote success probability of $X^{n}_{t}$. Hence,
\begin{eqnarray}   \label{binmod}
	X^{n}_{t} & \sim & \mbox{Binomial}(n_{0},p^{n}_{t}),
\end{eqnarray}
\begin{eqnarray*}
	p_{t}^{n} & = & \left\{ \begin{array}{ll}  p_{1} \neq 0.001, & t \geq \nu, n \in \mathcal{N},\\
		p_{0}=0.001, & \mbox{otherwises}. \end{array} \right.
\end{eqnarray*}

Let $k = t-s+1, S^{n}_{st} = \sum\limits^{t}_{i=s}X^{n}_{i} \sim \mbox{Binomial}(kn_{0},p_{t}^{n}).$ P-values are calculated based on $p^{n}_{st} = 2 \min (\varphi^{n}_{st}, 1-\varphi^{n}_{st})$, where $\varphi^{n}_{st} \sim$ Uniform$(P( S < S^{n}_{st}),P(S \leq S^{n}_{st}))$, and $S \sim \mbox{Binomial}(kn_{0},p_{0}).$

We select $n_{0} = 5, p_{0} =0.001, p_{1}= 0.05, \lambda_{1}=1, \lambda_{2}= \sqrt{\frac{\log \gamma}{\log \log \gamma}} $. Let $k = t-s+1$ and $S^{n}_{st}=\sum\limits^{t}_{i=s}X^{n}_{i} \sim \mbox{Binomial}(kn_{0},p_{1})$ if $n \in \mathcal{N}$ and $\nu = 1$. p-values are computed and SL stopping rule is applied.





For the two non-normal models and examples in $\ref{Poi_eg}$ and $\ref{Bin_eg}$ 
we all consider simulation settings $N=100$, $\# \mathcal{N}$ ranging from $1$ to $100, \nu = 1, ARL = 5000$ and $\mathcal{K} = \{1,\cdots, 200\}$. The detection thresholds are in Table $\ref{table:25}$, the detection delays in Table $\ref{table:26}$.

As the results in Table $\ref{table:26}$ show, SL stopping rule is capable to detect online change-points on multiple data streams for non-normal models. 

\begin{table}[ht]  
	\centering
	\caption{Thresholds $C_{\gamma}$ for stopping rules calibrated to $ARL=5000$.}
	\begin{tabular}{l|ll}
		\hline
		& \multicolumn{2}{c}{$N=100$} \\
		\hline
		$Test$ & $C_{\gamma}$ & $ARL$ \\
		\hline
		$SL($Poi,$\lambda_{1}=1,\lambda_{2}=\sqrt{\frac{\log \gamma}{\log\log \gamma}} )$ & 9.1 & 4865  \\
		$SL($Bin,$\lambda_{1}=1,\lambda_{2}=\sqrt{\frac{\log \gamma}{\log\log \gamma}})$ & 9.1 & 5072  \\
		\hline
	\end{tabular}
	\label{table:25}
\end{table}

\begin{table}[ht]    
	\centering
	\caption{Detection delays when $\# \mathcal{N}$ (out of $N = 100$) data streams undergo Poisson and binomial distribution changes.}
	\begin{tabular}{l|lllllll}
		\hline
		& \multicolumn{7}{|c}{$\# \mathcal{N}$}\\
		\hline
		$Test$ & 1 & 3 & 5 & 10 & 30 & 50 & 100\\
		\hline
		$SL($Poi,$\lambda_{1}=1,\lambda_{2}=1.99)$ & 27.6 & 12.7 & 8.8 & 5.3 & 2.3 & 1.5 & 1.0\\
		$SL($Bin,$\lambda_{1}=1,\lambda_{2}=1.99)$ & 23.6 & 11.1 & 7.6 & 4.5 & 1.9 & 1.3 & 1.0\\		
		\hline
	\end{tabular}
	\label{table:26}
\end{table}

\clearpage

\section{Disclosure statement}

The author reports there are no competing interests to declare.

\clearpage

\begin{appendix}
	
\section{Proof of Theorems \ref{thm21}}
\label{Proof21}
We summarize below the definitions
of the probability measures used in the proofs of Theorems $\ref{thm21}-\ref{thm22}$ in this and the next section:\\
1.$P_{s}(E_{s})$: This is the probability measure (expectation) under which an arbitrarily
chosen sequence has probability $(1-\epsilon)$ that all observations are ($i.i.d.$)
$N(0,1)$, and probability $\epsilon$ that observations are $N(0,1)$ before time $s$, $N(\Delta, 1)$ at and after time $s$.
In particular, if\\
$\hspace{1cm}$(a) $s=\infty$, then with probability $1$ all observations are $N(0,1)$.\\
$\hspace{1cm}$(b) $s=1$, then an arbitrarily chosen sequence has probability $(1-\epsilon)$ that all observations are $N(0,1)$, and probability $\epsilon$ that all observations are $N(\Delta,1)$.\\
2. $P(E)$: This is the probability measure (expectation) under which $Y,Y_{1},Y_{2},\cdots$,
are $i.i.d. N(0,1)$ random variables.\\

{\sc Proof of Theorem} $\ref{thm21}(a)$. \

Please refer to Section 6 in $\cite{Chan17}$.

\bigskip
{\sc Proof of Theorem} $\ref{thm21}(b)$.\\

Let $E_{\nu, \mathcal{N}}$ denote expectation with respect to $X_{t}^{n} \sim N(\mu_{t}^{n},1),$ with $\mu_{t}^{n}$ satisfying (\ref{munor}). For a given stopping rule $T$, define
\begin{equation}
	D_{N}^{(V)}(T) = \sup_{1 \leq \nu < \infty} \Big[\max_{\mathcal{N}:\# \mathcal{N}=V } E_{\nu, \mathcal{N}}(T-\nu+1|T \geq \nu)\Big].
\end{equation}
Let 
\begin{equation} \label{1b1}
	k = \left \lfloor{(1-\epsilon)2\Delta^{-2}V^{-1}N^{\zeta}}\right \rfloor, 
\end{equation}
with some $\epsilon >0$ and $V=o(\frac{N^{\zeta}}{\log N})$.
By Lemma 1 in $\cite{Chan17}$, we can find $s \geq 1$ such that
\begin{equation} \label{1b2}
	P_{\infty} (T \geq s+k|T\geq s ) \geq 1-k/\gamma. 
\end{equation}
Let $t = s+k-1,$ and consider the test, conditioned on $T\geq s,$ of
\begin{eqnarray*}  
	H_{0} & : & X_{u}^{n} \sim N(0,1) \mbox{ for } 1\leq n \leq N, 1 \leq u \leq t,\\
	\mbox{vs } H_{s,V} & : & X_{u}^{n} \sim N(\Delta \mathbf{1}_{\{u \geq s, n \in \mathcal{N}\} },1),
\end{eqnarray*}
\smallskip
for $1\leq n \leq N,$ $1 \leq u \leq t$, with $\mathcal{N}$ a random subset of $\{1, \cdots,N\}$ of size $V$.

By $(\ref{1b2}),$ the test rejecting $H_{0}$ when $T<s+k$ has Type I error probability not exceeding $k/\gamma$.

Let $A_{j} = \{\mathcal{N}:\# \mathcal{N} = j\}$. At time $t,$ the (conditional) likelihood ratio between $H_{s,V}$ and $H_{0}$ is $L_{V},$ where
\begin{equation*} 
	L_{j} = \binom{N}{j}^{-1}\sum_{\mathcal{N} \in A_{j}} \Big(\prod_{n \in \mathcal{N}}e^{Z^{n}_{st}\Delta\sqrt{k}-k\Delta^{2}/2} \Big). 
\end{equation*}
Let $P_{s,V}(E_{s,V})$ denote probability (expectation) with respect to $H_{s,V}$.
\bigskip
Let $\epsilon_{1}=2N^{-\beta}$ and
\begin{eqnarray}  \label{1b3}
	L(\epsilon_{1})&=&\prod_{n=1}^{N}(1-\epsilon_{1}+\epsilon_{1}e^{Z^{n}_{st}\Delta \sqrt{k}-k\Delta^{2}/2})\\ \nonumber
	&\Big[=&\sum_{j=0}^{N}(1-\epsilon_{1})^{N-j}\epsilon_{1}^{j} \binom{N}{j} L_{j} \Big].
\end{eqnarray}
Since $Z^{n}_{st} \sim N(\Delta \sqrt{k},1)$ if $n \in \mathcal{N}$ and $Z^{n}_{st} \sim N(0,1)$ if $n \notin \mathcal{N}$, it follows that
\begin{eqnarray*} 
	Ee^{Z^{n}_{st}\Delta \sqrt{k}-k\Delta^{2}/2} = \left\{ \begin{array}{ll} e^{k \Delta^{2}} & \mbox{if } n \in \mathcal{N}, \\
		1 & \mbox{if } n \notin \mathcal{N}. \end{array} \right.
\end{eqnarray*}
Therefore, by $(\ref{1b3}),$
\begin{equation}   \label{1b4}
	E_{s,V}L(\epsilon_{1})=(1-\epsilon_{1}+\epsilon_{1}e^{k\Delta^{2}})^{V}.
\end{equation}
By the monotonicity $E_{s,V}L_{1} \leq \cdots \leq E_{s,V}L_{N},$ and by $P(W\geq V) \rightarrow 1$ for $W \sim$ Binomial$(N, \epsilon_{1}),$ it follows from $(\ref{1b3})$ that
\begin{equation}   \label{1b5}
	E_{s,V}L(\epsilon_{1}) \geq P(W \geq V )E_{s,V}L_{V}=[1+o(1)]E_{s,V}L_{V}.
\end{equation}

For $\beta > 1-\zeta, k = \left \lfloor{(1-\epsilon)2\Delta^{-2}V^{-1}N^{\zeta}}\right \rfloor,$ with some $\epsilon >  0$. We show that
\begin{equation}    \label{1b6}
	(1-\epsilon_{1}+\epsilon_{1}e^{k\Delta^{2}})^{V} \leq (2\epsilon_{1}e^{k\Delta^{2}})^{V} \leq e^{k\Delta^{2}V}.
\end{equation}
We can choose $\epsilon>0$ such that 
\begin{equation}    \label{1b7} 
	e^{k\Delta^{2}V} = o(\exp(2N^{\zeta}/3)).
\end{equation}
By (\ref{1b6}) and (\ref{1b7}),
\begin{equation}   \label{1b8}
	(1-\epsilon_{1}+\epsilon_{1}e^{k\Delta^{2}})^{V}=o(\exp(2N^{\zeta}/3)).
\end{equation}
By $(\ref{1b4}),(\ref{1b5}), (\ref{1b8})$ and Markov's inequality, \begin{equation}   \label{1b9}
	P_{s,V}(L_{V} \geq J) \rightarrow 0, J = \exp(2N^{\zeta}/3).
\end{equation}
Let $c_{\gamma}$ be such that $P_{s,V}(J \geq L_{V} \geq c_{\gamma}) = \exp(-N^{\zeta}/4)$. It follows from $(\ref{1b9})$ that
\begin{equation}   \label{1b10}
	P_{s,V}(L_{V}\geq c_{\gamma} )(=P_{s,V}(L_{V}\geq c_{\gamma}|T \geq s )) \rightarrow 0,
\end{equation}
and that for $N$ large,
\begin{eqnarray} \label{1b11}
	P_{\infty}(L_{V} \geq c_{\gamma})(=P_{\infty}(L_{V} \geq c_{\gamma}|T \geq s)) 
	& \geq & P_{\infty}(J \geq L_{V} \geq c_{\gamma}) \\ \nonumber
	& = & E_{s,V}(L_{V}^{-1} \mathbf{1}_{\{J \geq L_{V} \geq c_{\gamma}\}}) \\ \nonumber
	& \geq & J^{-1}\exp(-N^{\zeta}/4) \geq k/\gamma.
\end{eqnarray}
By $(\ref{1b2}), (\ref{1b11})$ and the Neyman-Pearson lemma, the test rejecting $H_{0}$ when $L_{V} \geq c_{\gamma}$ is at least as powerful as the one based on $T,$ that is, 
\begin{equation}   \label{1b12}
	P_{s,V}(T \geq s+k|T \geq s) \geq P_{s,V}(L_{V} < c_{\gamma} ).
\end{equation}
It follows from $(\ref{1b10})$ and $(\ref{1b12})$ that
\begin{eqnarray*}  
	D_{N}^{(V)}(T) & \geq & E_{s,V}(T-s+1|T \geq s)\\
	& \geq & kP_{s,V}(T \geq s+k|T\geq s ) \\
	& = & k[1+o(1)],
\end{eqnarray*}
and the proof of Theorem $\ref{thm21}$(b) is complete.

\bigskip{}
\section{Proof of Theorems \ref{thm22}}    \label{Proof22}
\begin{lem} \label{lem24}
	Consider stopping rule $T_{SL}$, with window sizes in $(\ref{win})$. If the stopping rule threshold $C_{\gamma} = \log(4\gamma^{2}+2\gamma) \sim 2\log(\gamma),$ then $E_{\infty}T_{SL} \geq \gamma.$
\end{lem}

{\sc Proof}. 
By Markov inequality. it suffices to show that
\begin{equation}   \label{2_1}
	P_{\infty}(T_{SL} < 2\gamma) \leq \frac{1}{2}.
\end{equation}

Let $k=t-s+1, S^{n}_{st}=\sum\limits_{i=s}^{t}X_{i}^{n}, Z^{n}_{st}=\frac{S^{n}_{st}}{\sqrt{k}},$ and $ p^{n}_{st}=2\Phi(-|Z^{n}_{st}|)$.

Since $\ell(\mathbf{p_{st}})=\sum\limits_{n=1}^{N}\ell(p^{n}_{st})$ is a log likelihood ratio against $p^{n}_{st} \stackrel{i.i.d.}{\sim}$ Uniform$(0,1),$ it follows from a change of measure argument that 
\begin{equation}
	P_{\infty}(\ell(\mathbf{p_{st}}) \geq C_{\gamma}) \leq e^{-C_{\gamma}} = (4\gamma^{2}+2\gamma)^{-1}.
\end{equation}
By Bonferroni's inequality, 
\begin{eqnarray}
	& & P_{\infty}(T_{SL}<2\gamma) \\ \nonumber &\leq &  \sum_{(s,t):1 \leq s \leq t < 2\gamma} P_{\infty}(\ell(\mathbf{p_{st}}) \geq C_{\gamma}) \\ \nonumber
	& \leq &  {\left \lfloor{2\gamma+1}\right \rfloor \choose 2}\frac{1}{2\gamma(1+2\gamma)} \\ \nonumber
	&\leq & \frac{1}{2}, 
\end{eqnarray}
and $(\ref{2_1})$ follows.

\begin{lem} ($\cite{Hu2022}$) \label{lem25}
	Let $\bq = (q^1, \ldots, q^N)$,
	with $q^n \stackrel{i.i.d.}{\sim}$ {\rm Uniform(0,1)}.
	For fixed $\lambda_1 \geq 0$ and $\delta > 0$,
	$$\sup_{\delta \leq \lambda_2 \leq \sqrt{N}} P ( \ell(\bq) \leq -C \lambda_2^2 ) \rightarrow 0 \mbox{ as } C \rightarrow \infty \mbox{ and } N \rightarrow \infty.
	$$
\end{lem}

Assume $(\ref{gam}),(\ref{ep})$ and let the threshold $C_{\gamma}[\leq \log(4\gamma^{2}+2\gamma)$ by Lemma $\ref{lem24}$] be such that $E_{\infty}T_{SL}=\gamma.$ Let $\eta = \min\limits_{V\in J_{N}} P_{1}(\sum\limits_{n=1}^{N}\ell(p^{n}_{1k})\geq C_{\gamma}|\# \mathcal{N}=V),$ where
\begin{eqnarray} \label{JN}
	J_{N} = \left\{ \begin{array}{ll} \{V:V=o(\frac{N^{\zeta}}{\log N})  \}, & \mbox{if }\beta > 1-\zeta , \\
		\{V: V \geq \frac{1}{2}N\epsilon\}, & \mbox{if }\beta < 1-\zeta.  \end{array} \right.
\end{eqnarray}
In $\cite{Chan17}$, we know that by applying large deviation theory, 
\begin{equation}  \label{lar_dev}
	P_{1}(\#\mathcal{N} < \frac{1}{2}N\epsilon) \leq \exp(-\frac{1}{8}N\epsilon).
\end{equation}
We show in various cases below that $\eta \rightarrow 1.$ When
\begin{eqnarray*}
	k = \left\{ \begin{array}{ll} 1, & \mbox{if }\beta < \frac{1-\zeta}{2}, \\
		\left \lfloor{2(1+\epsilon)\Delta^{-2}\rho_{Z}(\beta, \zeta)\log N}\right \rfloor, & \mbox{if }\frac{1-\zeta}{2} < \beta < 1-\zeta, \\ \left \lfloor{2(1+\epsilon)\Delta^{-2}V^{-1} N^{\zeta}}\right \rfloor, & \mbox{if } \beta > 1-\zeta, \end{array} \right.
\end{eqnarray*}
for some $\epsilon>0$.
For $j \geq 1 $ and $V \in J_{N},P_{1}(T_{SL} \geq jk+1|\#\mathcal{N}=V ) \leq (1-\eta)^{j}.$ Hence, if $\beta < 1-\zeta,$ then $N\epsilon \sim N^{1-\beta} >> N^{\zeta}$ and it follows from (\ref{lar_dev}) that
\begin{equation}
	D_{N}(T_{SL}) \leq k\sum_{j=0}^{\infty}(1-\eta)^{j}+\gamma P_{1}(\# \mathcal{N} \notin J_{N}) \sim k.
\end{equation}
Similarly, if $V=o(\frac{N^{\zeta}}{\log N})$ for $\beta > 1-\zeta$, then

\begin{equation}
	D_{N}^{(V)}(T_{SL}) \leq [1+o(1)]k.
\end{equation}

For notational simplicity, we shall let $C$ denote a generic positive constant. Let $\textbf{q} = (q^{1}, \cdots, q^{N}),$ with $q^{n} \stackrel{i.i.d.}{\sim}$ Uniform$(0,1).$ Let $\nu$ be the change-point. For fixed $\lambda_{1} \geq 0$ and $\delta>0,$
by Lemma $\ref{lem25}$,
\begin{equation} \label{eta21}
	\sup_{\delta \leq \lambda_2 \leq \sqrt{N}} P ( \ell(\bq) \leq -C \lambda_2^2 ) \rightarrow 0 \mbox{ as } C \rightarrow \infty \mbox{ and } N \rightarrow \infty.
\end{equation}
\begin{equation} \label{eta22}
	P_{1} \Big(\sum_{n:n \notin \mathcal{N} } \ell(q^{n}) \geq -C\log \gamma \Big) \rightarrow 1 \mbox{ for any } C>0.
\end{equation}

We show $\eta \rightarrow 1$ in various cases below:\\

Case 1: $k=1$ for $\beta < \frac{1-\zeta}{2}$. 
When $n \in \mathcal{N}, $
\begin{equation} \label{E1}
	|E_{1}Z^{n}|=|\Delta|,
\end{equation}
with $Z^{n}=Z^{n}_{st}$. 

Let $\Gamma = \{n: n \in \mathcal{N}, |Z^{n}| \geq |\Delta| \}.$ Let $\# \Gamma \sim $Binomial$(N,r_{N})$, where
\begin{equation} \label{r1}
	r_{N}  =  N^{-\beta}\Phi(-|\Delta|+|\Delta|) = \frac{1}{2}N^{-\beta}.
\end{equation}
If $Q^{n}=1$ and $|Z^{n}| \geq z_{0}$ for $z_{0} \leq |E_{1}Z^{n}|,$
\begin{eqnarray} \label{pn1}
	p^n & = &  2\int_{|Z^n|}^{\infty} \frac{1}{\sqrt{2 \pi}} e^{-y^{2}/2} dy \\ \nonumber
	& = &  2\int_{|Z^n|-z_0}^{\infty} \frac{1}{\sqrt{2 \pi}} e^{-(z+z_0)^2/2} dz \leq  2e^{-z_0^2/2} q^n = o(q^{n}),
\end{eqnarray}
with $p^{n}=p^{n}_{st}.$ By $(\ref{E1}),$ for $n \in \Gamma$,
\begin{eqnarray} \label{pn}
	p^{n} \leq 2e^{-\frac{\Delta^{2}}{2}}.
\end{eqnarray}
Let $2e^{-\frac{\Delta^{2}}{2}}=a$ and by $(\ref{pn1})$, $\lambda_2 \sim \sqrt{\frac{N^{\zeta}}{\zeta \log N}}$ and $\log(1+x) \sim x$ as $x \rightarrow 0$,
\begin{equation} \label{lpq1}
	\sum_{n \in \Gamma} [\ell(p^n)-\ell(q^n)] \geq [1+o(1)] \frac{\lambda_2 (\# \Gamma)}{\sqrt{aN \log N}} \sim \frac{\sqrt{N^{\zeta}} (\# \Gamma)}{\sqrt{Na\zeta} \log N}.
\end{equation}
We also know that 
\begin{equation} \label{cGam}
	C_{\gamma} = \log(4\gamma^{2}+2\gamma) \sim 2\log(\gamma) \sim 2N^{\zeta}.
\end{equation}
Since $\ell(p^n) \geq \ell(q^n)$ for all $n$, $\eta \rightarrow 1$ follows from (\ref{eta22}), (\ref{r1}), (\ref{lpq1}) and (\ref{cGam}).

\bigskip
Case 2: $k = \left \lfloor{2(1+\epsilon)\Delta^{-2}\rho_{Z}(\beta, \zeta)\log N}\right \rfloor$ for $\frac{1-\zeta}{2} < \beta < \frac{3(1-\zeta)}{4}$.
When $n \in \mathcal{N}, $
\begin{eqnarray} \label{E22}
	|E_{1}Z^{n}| &=&|\Delta|\sqrt{k} \geq \sqrt{2\nu \log N},\\ \nonumber
	\nu &=& (1+\epsilon)^{\frac{1}{2}}\rho_{Z}(\beta,\zeta),
\end{eqnarray}
with $Z^{n}=Z^{n}_{st}$. 

Let $\Gamma = \{n: n \in \mathcal{N}, |Z^{n}| \geq 2\sqrt{(2\beta-1+\zeta) \log N} \}.$ Let $\delta$ be such that $(0<)2\delta = \beta - \frac{1-\zeta}{2}-(2\sqrt{\beta-\frac{1-\zeta}{2}}-\sqrt{\nu})^{2}$. Let $\# \Gamma \sim$ Binomial$(N,r_{N}),$ where
\begin{eqnarray}  \label{Gam2}
	& & r_{N} \\ \nonumber
	& \geq & N^{-\beta}\Phi \Big( -2\sqrt{(2\beta-1+\zeta)\log N}+\sqrt{2 \nu \log N} \Big) \\ \nonumber
	& \geq & N^{-\beta-(2\sqrt{\beta-\frac{1-\zeta}{2}}-\sqrt{\nu})^2-\delta} = N^{\frac{1}{2}-2\beta-\frac{\zeta}{2}+\delta}.
\end{eqnarray}


Let $p^{n}$ and $q^{n}$ be defined in this section before. By (\ref{pn1}), for $n \in \Gamma$,
\begin{equation} \label{pn22}
	p^{n} \leq 2N^{-4\beta+2-2\zeta} \mbox{ and } p^{n} = o(q^n).
\end{equation}
By (\ref{pn22}), $\lambda_2 \sim \sqrt{\frac{N^{\zeta}}{\zeta \log N}}$ and $\log(1+x) \sim x$ as $x \rightarrow 0$,
\begin{equation} \label{lpq22}
	\sum_{n \in \Gamma} [\ell(p^n)-\ell(q^n)] \geq [1+o(1)] \frac{\lambda_2 N^{2\beta-1+\zeta}(\# \Gamma)}{\sqrt{2N \log N}} \sim \frac{N^{2\beta-\frac{3}{2}(1-\zeta)}(\# \Gamma)}{\sqrt{2\zeta} \log N}.
\end{equation}
Since $\ell(p^n) \geq \ell(q^n)$ for all $n$, $\eta \rightarrow 1$ follows from (\ref{eta22}), (\ref{cGam}),  (\ref{Gam2}) and (\ref{lpq22}).

\bigskip
Case 3: $k = \left \lfloor{2(1+\epsilon)\Delta^{-2}\rho_{Z}(\beta, \zeta)\log N}\right \rfloor$ for $\frac{3(1-\zeta)}{4} < \beta < 1-\zeta$.
When $n \in \mathcal{N}, $
\begin{eqnarray} \label{E3}
	|E_{1}Z^{n}|&=&|\Delta|\sqrt{k} \geq \sqrt{2\nu \log N},\\ \nonumber
	\nu &=& (1+\epsilon)^{\frac{1}{2}}\rho_{Z}(\beta,\zeta),
\end{eqnarray}
with $Z^{n}=Z^{n}_{st}$. 

Let $\Gamma = \{n: n \in \mathcal{N}, |Z^{n}| \geq \sqrt{2(1-\zeta) \log N} \}.$ Let $\delta$ be such that $(0<)2\delta = 1-\zeta -\beta-(\sqrt{1-\zeta}-\sqrt{\nu})^{2}$. Let $\# \Gamma \sim $Binomial$(N,r_{N}),$ where
\begin{eqnarray}  \label{Gam3}
	& & r_{N}  \\ \nonumber
	& \geq & N^{-\beta}\Phi \Big( -\sqrt{2(1-\zeta) \log N}+\sqrt{2 \nu \log N} \Big) \\ \nonumber
	& \geq & N^{-\beta-(\sqrt{1-\zeta}-\sqrt{\nu})^2-\delta} = N^{\zeta+\delta-1}.
\end{eqnarray}


Let $p^{n}$ and $q^{n}$ be defined in this section before. By (\ref{pn1}), for $n \in \Gamma$,
\begin{equation} \label{pn23}
	p^{n} \leq 2N^{\zeta-1} \mbox{ and } p^{n} = o(q^n).
\end{equation}
By (\ref{pn23}), $\lambda_2 \sim \sqrt{\frac{N^{\zeta}}{\zeta \log N}}$ and $\log(1+x) \sim x$ as $x \rightarrow 0$,
\begin{equation} \label{lpq23}
	\sum_{n \in \Gamma} [\ell(p^n)-\ell(q^n)] \geq [1+o(1)] \frac{\lambda_2 (\# \Gamma)}{\sqrt{2N^{\zeta} \log N}} \sim \frac{\# \Gamma}{\sqrt{2\zeta} \log N}.
\end{equation}
Since $\ell(p^n) \geq \ell(q^n)$ for all $n$, $\eta \rightarrow 1$ follows from (\ref{eta22}), (\ref{cGam}),  (\ref{Gam3}) and (\ref{lpq23}).

\bigskip
Case 4: $k = \left \lfloor{2(1+\epsilon)\Delta^{-2} V^{-1} N^{\zeta}}\right \rfloor$ for $ \beta > 1-\zeta$ and $V=o(\frac{N^{\zeta}}{\log N})$. When $n \in \mathcal{N}, $
\begin{equation} \label{E4}
	|E_{1}Z^{n}|=|\Delta|\sqrt{k} \geq \sqrt{2(1+\epsilon)^{\frac{1}{2}}V^{-1} N^{\zeta}},
\end{equation}
with $Z^{n}=Z^{n}_{st}$.\\

Let $\Gamma = \{n: n \in \mathcal{N}, |Z^{n}| \geq \sqrt{2(1+\epsilon)^{\frac{1}{3}} V^{-1}N^{\zeta}} \}.$ Since $\# \mathcal{N} = V,$ by ($\ref{E4})$ and the law of large number,
\begin{eqnarray} \label{Gam4}
	P_{1}(\# \Gamma \geq (1+\epsilon)^{-\frac{1}{4}}V) \rightarrow 1.
\end{eqnarray}

Let $p^{n}$ and $q^{n}$ be defined in this section before. By (\ref{E4}), for $n \in \Gamma$,
\begin{equation} \label{pn4}
	p^{n} \leq 2\exp(-(1+\epsilon)^{\frac{1}{2}}V^{-1} N^{\zeta}).
\end{equation}
Since $V^{-1} N^{\zeta}$ is large compared to $\log N,$ by (\ref{pn4}),
\begin{equation} 
	\ell(p^{n}) \geq (1+\epsilon)^{\frac{1}{3}} V^{-1} N^{\zeta} \mbox{  for  } n \in \Gamma.
\end{equation}
Moreover $\ell(p) \geq -1$ for $N$ large and therefore,
\begin{equation} \label{lp}
	\sum_{n: n \in \mathcal{N}} \ell(p^n) \geq (\# \Gamma)(1+\epsilon)^{\frac{1}{3}}V^{-1} N^{\zeta} - V.
\end{equation}
Hence, $\eta \rightarrow 1$ follows from (\ref{eta22}), (\ref{cGam}),  (\ref{Gam4}), (\ref{lp}), and $V^{-1}N^{\zeta} \rightarrow \infty$.

\end{appendix}

\end{document}